\documentclass{amsart}
\usepackage{amsmath,amsthm,amssymb,amscd,enumitem,mathtools}
\usepackage{amsfonts}
\usepackage{rotating}
\usepackage{euscript}
\usepackage{pst-node}
\usepackage{epsfig,verbatim}
\usepackage{euscript}
\mathtoolsset{showonlyrefs}
\def\be{\begin{equation}}
\def\ee{\end{equation}}

\def\C{{\mathbb C}} 
\def\f{\EuScript}
 
\def\P{{\mathbb P}}

\def\t{\widetilde}

\def\ord{{\rm ord}}

\def\phi{{\varphi}}

\def\deg{{\rm deg\,}}

\def\GCD{{\rm GCD }}

\def\bp{\begin{proposition}}
\def\ep{\end{proposition}}

\def\bt{\begin{theorem}}
\def\et{\end{theorem}}
\def\br{\begin{remark}}
\def\er{\end{remark}}
\def\be{\begin{equation}}
\def\bee{\begin{equation*}}
\def\l{\label}

\def\c{\mathcal}
\def\ee{\end{equation}}
\def\eee{\end{equation*}}
\def\bl{\begin{lemma}}
\def\el{\end{lemma}}
\def\bc{\begin{corollary}}
\def\ec{\end{corollary}}
\def\pr{\noindent{\it Proof. }}

\def\bd{\begin{definition}}
\def\ed{\end{definition}}
\def\t{\widetilde}

\newtheorem{theorem}{Theorem}[section]
\newtheorem{lemma}[theorem]{Lemma}
\newtheorem{definition}[theorem]{Definition}
\newtheorem{corollary}[theorem]{Corollary}
\newtheorem{proposition}[theorem]{Proposition}
\newtheorem{problem}[theorem]{Problem}

\theoremstyle{definition}

\theoremstyle{definition}
\newtheorem{remark}[theorem]{Remark}
\def\bpr{\begin{problem}}
\def\epr{\end{problem}}

\begin{document}
\title{On algebraic dependencies between Poincar\'e functions}
\author{Fedor Pakovich}
\thanks{
This research was supported by ISF Grant No. 1092/22}
\address{Department of Mathematics, Ben Gurion University of the Negev, Israel}
\email{
pakovich@math.bgu.ac.il}

\begin{abstract} 
Let $A$ be a rational function of one complex variable of degree at least two, and $z_0$ its repelling fixed point
with the multiplier $\lambda.$  A Poincar\'e function associated with $z_0$ is a function $\f P_{A,z_0,\lambda}$  meromorphic on $\C$ such that  
$\f P_{A,z_0,\lambda}(0)=z_0$, $\f P_{A,z_0,\lambda}'(0)\neq 0,$   and $\f P_{A,z_0,\lambda}(\lambda z)=A\circ \f P_{A,z_0,\lambda}(z).$
In this paper, we study the following problem:  given 
Poincar\'e functions $\f P_{A_1,z_1,\lambda_1}$ and $\f P_{A_2,z_2,\lambda_2}$,  
find out if there is an algebraic relation $f(\f P_{A_1,z_1,\lambda_1},\f P_{A_2,z_2,\lambda_2})=0$
between them and, if such a relation exists, describe the corresponding algebraic curve $f(x,y)=0.$ We provide a solution, which can be viewed as a refinement of the classical theorem of Ritt about commuting rational functions. We also reprove and extend previous results concerning algebraic dependencies between B\"ottcher functions.

\end{abstract} 
\maketitle

\section{Introduction}
Let $A$ be a rational function of one complex variable of degree at least two, and $z_0$ its repelling fixed point
with the multiplier $\lambda.$ 
We recall that a {\it Poincar\'e function} $\f P_{A,z_0,\lambda}$ associated with $z_0$ is a function meromorphic on $\C$ such that 
$\f P_{A,z_0,\lambda}(0)=z_0,$ $\f P_{A,z_0,\lambda}'(0)\neq 0,$  
and the diagram    
\be \l{xor}
\begin{CD} 
\C @>\lambda z>> \C \\ 
@V \f P_{A,z_0,\lambda} VV  @VV \f P_{A,z_0,\lambda} V\\ 
\f \C\P^1 @>A>> \f\C\P^1\
\end{CD}
\ee
commutes. 
The Poincar\'e function  exists and is defined  up to 
the transformation of argument $z\rightarrow cz,$ where $c\in \C^*$ (see e. g. \cite{mil}). 
In particular, it is defined in a unique way if to assume that 
$\f P_{A,z_0,\lambda}^{\prime}(0)=1$. Such Poincar\'e functions are called  normalized.  
In this paper, we will consider non-normalized Poincar\'e functions, so the explicit meaning of the notation $\f P_{A,z_0,\lambda}$ is following: $\f P_{A,z_0,\lambda}$ is 
{\it some}  meromorphic function satisfying the above conditions. 
We say that a rational function  $A$ is {\it special} if it is either a Latt\`es map, or it is conjugate to $z^{\pm n}$ or $\pm T_n.$  Poincar\'e functions associated with special functions can be described in terms of classical functions. Moreover, by the result of Ritt \cite{r0}, these functions are the only  Poincar\'e functions that are periodic. 

In this paper, we study the following problem.
Let  $A_1$, $A_2$ be non-special rational functions of degree at least two with  repelling fixed points $z_1$, $z_2$, and $ \f P_{A_1,z_1,\lambda_1}$, $\f P_{A_2,z_2,\lambda_2}$ corresponding Poincar\'e functions. Under what conditions    there exists an algebraic curve $f(x,y)=0$ such that  \be \l{xre}f(\f P_{A_1,z_1,\lambda_1},\f P_{A_2,z_2,\lambda_2})=0\ee and, if such a curve exists, how it can be described? 
The simplest example of relation \eqref{xre} is just the equality 
\be \l{sim} \f P_{A_1,z_0,\lambda_1}=\f P_{A_2,z_0,\lambda_2},\ee which is known to 
have strong dynamical consequences. 
Specifically, equality \eqref{sim} implies easily that $A_1$ and $A_2$ commute. On the other hand, by the theorem of Ritt  (see \cite{r} and also \cite{e2}, \cite{pak}), every two non-special commuting rational functions of degree at least two have a common iterate. Thus, equality \eqref{sim} implies that 
\be \l{ite} 
A_1^{\circ l_1}=A_2^{\circ l_2}
\ee
for some integers $l_1,l_2\geq 1$. Moreover, the Ritt theorem essentially is equivalent to the statement that equality \eqref{sim} implies equality \eqref{ite}, since it was observed already by Fatou and Julia (\cite{f}, \cite{j})  that if two rational functions commute, then 
some of their iterates share a repelling fixed point and a corresponding Poincar\'e function.   

To our best knowledge, the problem of describing algebraic dependencies between Poincar\'e functions has never been considered in the literature. Nevertheless,  
the problem of describing  algebraic dependencies between {\it B\"ottcher functions}, similar in spirit, has been investigated in the papers \cite{bb1}, \cite{ng}. We recall that for a polynomial $P$ of degree $n$ a corresponding  B\"ottcher function $\f B_P$ is a Laurent series 
\be \l{boet1} \f B_P=a_{-1} z+a_0+\frac{a_1}{z}+\frac{a_2}{z^2}+\dots \in z\C [[1/z]], \ \ \ a_{-1}\neq 0,\ee that makes the diagram 
\be \l{boet2} 
\begin{CD} 
\C @>z^n>> \C \\ 
@V  \f B_A VV  @VV \f B_A V\\ 
\f \C\P^1 @>A>> \f\C\P^1\
\end{CD}
\ee 
commutative. In this notation, the result of Becker and Bergweiler \cite{bb1} (see also \cite{bb2}),  
 states that if $A_1$ and $A_2$ are polynomials of the same degree $d$, then  the function 
$\beta=\f B_{A_1}\circ \f B_{A_2}^{-1}$ is transcendental, unless either $\beta$ is linear, or $A_1$ and $A_2$ are special (notice that since a polynomial cannot be a Latt\`es map, a polynomial is special if and only if it is conjugate to $z^n$ or $\pm T_n$).  Since the equality 
\be \l{ses} f(\f B_{A_1}(z), \f B_{A_2}(z))=0\ee holds for some $f(x,y)\in \C[x,y]$ if and only if the function $\beta$ is algebraic, this result implies the absence of  algebraic dependencies of degree greater than one between $\f B_{A_1}(z)$ and $\f B_{A_2}(z)$ for non-special $A_1$ and $A_2$ of the same degree. 
 
Subsequently, it was proved by Nguyen in the paper \cite{ng}   that 
the equality 
\be \l{es} f(\f B_{A_1}(z^{d_1}), \f B_{A_2}(z^{d_2}))=0\ee 
holds  for  some integers $d_1,d_2\geq 1$ if and only if there exist  polynomials $X_1,X_2, B$  and integers $l_1, l_2\geq 1$  such that the diagram 
\be \l{hol+} 
\begin{CD} 
 (\C\P^1)^2 @> (B,B) >> (\C\P^1)^2
 \\ 
@V (X_1,X_2) VV @VV  (X_1,X_2)   V\\ 
(\C\P^1)^2 @> (A_1^{\circ l_1}, A_2^{\circ l_2}) >> (\C\P^1)^2 
\end{CD} 
\ee
commutes. Notice that although the result of Nguyen deals with the more general situation than  the result of Becker and Bergweiler,  the former does not formally imply the latter.

Let us recall that an algebraic curve $C: f(x,y)=0$ has genus zero if and only if it admits a parametrization $z\rightarrow (X_1(z),X_2(z))$ by rational functions $X_1,$ $X_2.$ Such a parametrization is called {\it generically one-to-one} if it is one-to-one except for finitely many points. By the L\"uroth theorem, this equivalent to say  that 
$X_1$ and $X_2$ generate the whole field of rational functions $\C(z)$. In this notation,    
our main result is the following analogue of the result of Nguyen.

\bt \l{m2} Let $A_1$, $A_2$ be  non-special rational functions of degree at least two,  $z_1$, $z_2$ their repelling fixed points with multipliers $\lambda_1,$ $\lambda_2$, 
and $\f P_{A_1,z_1,\lambda_1}$, $\f P_{A_2,z_2,\lambda_2}$ Poincar\'e functions. 
Assume that  $C: f(x,y)=0$ is an irreducible algebraic curve, and $d_1,$ $d_2$ are coprime positive integers such that the equality 
\be \l{xru} f\left(\f P_{A_1,z_1,\lambda_1}(z^{d_1}),\f P_{A_2,z_2,\lambda_2}(z^{d_2})\right)=0\ee holds. Then $C$ has genus zero. Furthermore, if $C: f(x,y)=0$ is an irreducible algebraic curve of genus zero with a generically one-to-one parametrization by rational functions 
$z\rightarrow (X_1(z),X_2(z))$, and $d_1,$ $d_2$ are coprime positive integers, then  equality \eqref{xru} holds for some  Poincar\'e functions $\f P_{A_1,z_1,\lambda_1}$, $\f P_{A_2,z_2,\lambda_2}$ if and only if  
 there exist positive integers $l_1,l_2,k$ and a rational function $B$ with a repelling fixed point  $z_0$  such that the diagram 
\be \l{hol} 
\begin{CD} 
 (\C\P^1)^2 @> (B,B) >> (\C\P^1)^2
 \\ 
@V (X_1,X_2) VV @VV  (X_1,X_2)   V\\ 
(\C\P^1)^2 @> (A_1^{\circ l_1}, A_2^{\circ l_2}) >> (\C\P^1)^2, 
\end{CD} 
\ee
commutes and the equalities 
\be \l{io} X_1(z_0)=z_1, \ \ \ \ X_2(z_0)=z_2,\ee  
\be \l{oi} \ord_{z_0}X_1=d_1k, \ \ \ \  \ord_{z_0}X_2=d_2k \ee
hold. 
\et

Notice that Theorem \ref{m2} can be considered as a refinement of the Ritt theorem.  
Indeed, equality \eqref{sim} is a particular case of the condition \eqref{xru}, where the curve $$f(x,y)=x-y=0$$ is parametrized by the functions $X_1=z,$ $X_2=z.$ Thus, in this case diagram \eqref{hol} reduces to equality \eqref{ite}. 
 More generally, considering the curve $x-R(y)=0,$ where $R$ is a rational function, we conclude that the equality 
$$\f P_{A_1,z_1,\lambda_1}=R\circ {\f P}_{A_2,z_2,\lambda_2}$$ implies that  
there exist $l_1,l_2\geq 1$ such that the diagram 
$$ 
\begin{CD} 
\C\P^1 @> A_2^{\circ l_2}>> \C\P^1 \\
@VV R V @VV  R V\\ 
\f \C\P^1 @>A_1^{\circ l_1}>> \f\C\P^1
\end{CD}
$$
 commutes.

Notice also that Theorem \ref{m2} implies the following handy criterion for the algebraic independence of Poincar\'e functions. 

\bc \l{cor1} Let $A_1$, $A_2$ be non-special rational functions of degrees $n_1\geq 2$, $n_2\geq 2$, 
 and  $z_1$, $z_2$ their repelling fixed points with multipliers $\lambda_1,$ $\lambda_2$. Then  
Poincar\'e functions $\f P_{A_1,z_1,\lambda_1}$, $\f P_{A_2,z_2,\lambda_2}$ 
 are algebraically independent, unless 
there exist positive integers $l_1,l_2$ and $l_1',l_2'$
such that $n_1^{l_1}=n_2^{l_2}$ and $\lambda_1^{l_1'}=\lambda_2^{l_2'}$. 
\ec

In addition to Theorem \ref{m2}, we prove the following more precise version of the theorem of Nguyen, which formally includes and generalizes 
the result of Becker and Bergweiler.

\bt \l{m3} Let $A_1$, $A_2$ be  non-special polynomials of degree at least two, and 
$\f B_{A_1}$, $\f P_{A_2}$ B\"ottcher functions.  Assume that  $C: f(x,y)=0$ is an irreducible algebraic curve, and $d_1,$ $d_2$ are coprime positive integers such that the equality 
\be \l{xru2} f\left(\f B_{A_1}(z^{d_1}),\f B_{A_2}(z^{d_2})\right)=0\ee holds. Then $C$ has the form $Y_1(x)-Y_2(y)=0$, where $Y_1,Y_2$  are  polynomials of coprime degrees, and can be parametrized by polynomials. 
Furthermore, if \linebreak $C: f(x,y)=0$ is an irreducible algebraic curve as above with a generically one-to-one parametrization by polynomials  
$z\rightarrow (X_1(z),X_2(z))$, and $d_1,$ $d_2$ are coprime positive integers, then  equality \eqref{xru2} holds for some   B\"ottcher functions $\f B_{A_1}$, $\f B_{A_2}$  if and only if  
 there exist positive integers $l_1,l_2$ and a polynomial $B$   such that the diagram 
\be \l{hol2} 
\begin{CD} 
 (\C\P^1)^2 @> (B,B) >> (\C\P^1)^2
 \\ 
@V (X_1,X_2) VV @VV  (X_1,X_2)   V\\ 
(\C\P^1)^2 @> (A_1^{\circ l_1}, A_2^{\circ l_2}) >> (\C\P^1)^2, 
\end{CD} 
\ee
commutes, and the equalities 
\be \l{oi2} \deg X_1=d_1, \ \ \ \  \deg X_2=d_2 \ee
hold. In particular, the equality 
$$f\left(\f B_{A_1}(z),\f B_{A_2}(z)\right)=0$$ implies that $C: f(x,y)=0$ has degree one and some iterates of 
$A_1$ and $A_2$ are conjugate. 
\et

Notice that the parameters $d_1,d_2$ appear in conclusions of both Theorem \ref{m2} and Theorem \ref{m3}. However, the condition \eqref{oi} is less restrictive than the condition \eqref{oi2}.   In particular, applying Theorem \ref{m3} for $d_1=d_2=1$ we conclude that 
algebraic dependencies between B\"otcher functions are essentially  trivial. On the other hand, 
algebraic dependencies between Poincar\'e  functions do exist (see Section \ref{sec3}). 

The approach of Nguyen to the study of algebraic dependencies \eqref{es} relies on the fact that 
such dependencies  give rise to {\it invariant algebraic curves} 
for endomorphisms 
 \be \l{ss0}
(A_1,A_2)
:\, ( \C\P^1)^2\rightarrow  ( \C\P^1)^2,\ee given by the formula
 \be \l{ss} (z_1,z_2)\rightarrow (A_1(z_1),A_2(z_2)),\ee 
where  $A_1$ and $A_2$ are polynomials. Say, for $A_1$ and $A_2$ of the  same degree $n$,  this can be seen immediately, since 
after substituting $z^n$ for $z$ into \eqref{es} we obtain the equality $$f(A_1\circ \f B_{A_1}(z^{d_1}), A_2\circ \f B_{A_2}(z^{d_2}))=0,$$ implying that 
$f(x,y)=0$ is $(A_1,A_2)$-invariant. 
Invariant curves for polynomial endomorphisms \eqref{ss0}   were  classified by Medvedev and Scanlon in the paper  \cite{ms}, and the proof of the theorem of
Nguyen relies crucially on this classification. 

Our approach the the study of algebraic dependencies \eqref{xre} is similar. However, instead of the paper \cite{ms} we use the results of the recent paper \cite{inv} providing a classification of invariant curves for   endomorphisms \eqref{ss} defined by arbitrary non-special {\it rational} functions $A_1,$ $A_2$. 
Notice that the paper \cite{ms} is based on  the Ritt theory of polynomial decompositions (\cite{rpc}), which  does not extend to rational functions. Accordingly, the approach of \cite{inv} is completely different and relies on the recent results 
\cite{semi}, \cite{rec},  \cite{dyna}, \cite{lattes},  \cite{fin} 
about {\it semiconjugate rational functions},  
which appear naturally in a variety of different contexts (see e. g. \cite{b},  \cite{e3}, \cite{i}, \cite{ms}, \cite{ng}, \cite{pj}, \cite{lattes}, \cite{aol}, \cite{inv}).

 This paper is organized as follows. In the second section, we review the notion of a {\it generalized Latt\`es map}, introduced in \cite{lattes}, and recall some results about semiconjugate rational functions and invariant curves proved in \cite{inv}. In the third section, we 
prove Theorem \ref{m2}.
We also show that for rational functions that are not generalized Latt\`es maps  equality \eqref{xru} under the condition $\GCD(d_1,d_2)=1$ implies the equality $d_1=d_2=1$
 (Theorem \ref{m1}). 
Finally, in the fourth section,  basing on results of the paper \cite{pj}, which complements some of results of \cite{ms},   we reconsider algebraic dependencies between B\"ottcher functions  and prove Theorem \ref{m3}.

\section{\l{s1} Generalized Latt\`es maps and invariant curves} 
\subsection{\l{ss1} Generalized Latt\`es maps and semiconjugacies} 
Let us recall that {\it a Riemann surface  orbifold} is a pair $\f O=(R,\nu)$ consisting of a Riemann surface $R$ and a ramification function $\nu:R\rightarrow \mathbb N$, 
which takes the value $\nu(z)=1$ except at isolated points. 
For an orbifold $\f O=(R,\nu)$, 
 the {\it  Euler characteristic} of $\f O$ is the number
\be \l{euler} \chi(\f O)=\chi(R)+\sum_{z\in R}\left(\frac{1}{\nu(z)}-1\right).\ee 
For   orbifolds $\f O_1=(R_1,\nu_1)$  and $\f O_2=(R_2,\nu_2)$, 
we write $ \f O_1\preceq \f O_2$ 
if $R_1=R_2$ and for any $z\in R_1$ the condition $\nu_1(z)\mid \nu_2(z)$ holds.

Let $\f O_1=(R_1,\nu_1)$  and $\f O_2=(R_2,\nu_2)$ be orbifolds, 
and let 
$f:\, R_1\rightarrow R_2$  be a holomorphic branched covering map. We say that $f:\,  \f O_1\rightarrow \f O_2$
is  a {\it covering map} 
{\it between orbifolds}
if for any $z\in R_1$ the equality 
\be \l{us} \nu_{2}(f(z))=\nu_{1}(z)\deg_zf\ee holds, where $\deg_zf$ is the local degree of $f$ at the point $z$.
If for any $z\in R_1$  
the weaker condition 
\be \l{uuss} \nu_{2}(f(z))\mid \nu_{1}(z)\deg_zf\ee
is satisfied,   we say that $f:\,  \f O_1\rightarrow \f O_2$ 
is a {\it holomorphic map}
 {\it between orbifolds}. 
If $f:\,  \f O_1\rightarrow \f O_2$ is a covering map between orbifolds with compact supports, then  the Riemann-Hurwitz 
formula implies that 
\be \l{rhor} \chi(\f O_1)= \chi(\f O_2) \deg f. \ee 
More generally,  if $f:\,  \f O_1\rightarrow \f O_2$ is a holomorphic map, then \be \l{iioopp} \chi(\f O_1)\leq \chi(\f O_2)\,\deg f, \ee and the equality is attained if and only if $f:\, \f O_1\rightarrow \f O_2$ is a covering map between orbifolds
(see \cite{semi}, Proposition 3.2).

Let $R_1$, $R_2$ be Riemann surfaces and 
$f:\, R_1\rightarrow R_2$ a holomorphic branched covering map. Assume that $R_2$ is provided with a ramification function $\nu_2$. In order to define a ramification function $\nu_1$ on $R_1$ so that $f$ would be a holomorphic map between orbifolds $\f O_1=(R_1,\nu_1)$ and $\f O_2=(R_2,\nu_2)$ 
we must satisfy condition \eqref{uuss}, and it is easy to see that
for any  $z\in R_1$ a minimum possible value for $\nu_1(z)$ is defined by 
the equality 
\be \l{rys} \nu_{2}(f(z))=\nu_ {1}(z)\GCD(\deg_zf, \nu_{2}(f(z)).\ee 
In case \eqref{rys} is satisfied for  any $z\in R_1$, we 
say that $f$ is {\it a  minimal holomorphic  map} 
between orbifolds 
$\f O_1=(R_1,\nu_1)$ and $\f O_2=(R_2,\nu_2)$.

We recall that  {\it a Latt\`es map} can be defined as a rational function $A$ such that $A:\f O\rightarrow \f O$ is a  covering self-map
for some orbifold $\f O$ on $\C\P^1$ (see \cite{mil2}, \cite{lattes}). Thus, $A$ is a Latt\`es map  if there exists an orbifold $\f O=(\C\P^1,\nu)$ 
such that  for any $z\in  \C\P^1$ the equality 
\be \l{uu0} \nu(A(z))=\nu(z)\deg_zA\ee holds. By  formula \eqref{rhor}, such $\f O$ necessarily satisfies $\chi(\f O)=0.$  
Following \cite{lattes}, we say that a rational function $A$ of degree at least two is 
a {\it generalized Latt\`es map} if there exists an orbifold $\f O=(\C\P^1,\nu)$, 
 distinct from the non-ramified sphere, 
such that  $A:\f O\rightarrow \f O$ is a minimal holomorphic self-map between orbifolds; that is, for any $z\in  \C\P^1$, the equality 
\be \l{uu} \nu(A(z))=\nu(z)\GCD(\deg_zA,\nu (A(z)))\ee holds.  By inequality \eqref{iioopp}, such $\f O$ satisfies $\chi(\f O)\geq 0$. 
Notice that any special rational function is a generalized Latt\`es map and that some iterate $A^{\circ l}$, $l\geq 1,$  of a rational function $A$ is a generalized Latt\`es map if and only if  
 $A$ is a generalized Latt\`es map (see \cite{inv}, Section 2.3). 

Generalized Latt\`es map are closely related to the problem of describing  semiconjugate rational functions, that is, rational functions that make the diagram 
\be \l{ssee} 
\begin{CD} 
 \C\P^1 @> B >> \C\P^1
 \\ 
@V X VV @VV  X  V\\ 
\C\P^1 @> A >> \C\P^1 
\end{CD} 
\ee 
commutative. For a general theory we refer the reader to the papers 
\cite{semi}, \cite{rec},  \cite{dyna}, \cite{lattes},  \cite{fin}. Below we need only the following two results, 
which are simplified reformulations  of Proposition 3.3 and Theorem 4.14 in \cite{inv}.  

The first result  states that if the function  $A$ in \eqref{ssee} is not   a generalized Lattès map, then \eqref{ssee} can be completed to a diagram of the very special form.

\bp \l{1} Let $A$ be a  rational function of degree at least two that is not a generalized Lattes map, and $X, B$ rational functions such that diagram \eqref{ssee} commutes.  Then there exists a rational function $Y$ such that the diagram
\be  \l{xx}
\begin{CD} 
 \C\P^1 @>B>> \C\P^1 \\ 
@V X  VV @VV X V\\ 
  \C\P^1 @> A >>  \C\P^1
\\ 
@V Y  VV @VV  Y V\\ 
  \C\P^1 @>B >>  \C\P^1 
\end{CD}
\ee
commutes, and the 
equalities 
\be \l{en2} Y  \circ X  =B^{\circ d} \ \ \ \ \ \ \ X \circ Y =A ^{\circ d},
\ee hold for some $d\geq 0$.
\ep

The second result relates an arbitrary non-special  rational function with some   
 rational function that is not a generalized Lattès map through the semiconjugacy relation.  

\bt \l{33} Let $A$ be a non-special rational function of degree at least two. Then there exist rational functions $\theta$ and $F$ such that $F$ is not a generalized Lattès map and the diagram 
\be \l{asa}
\begin{CD} 
 \C\P^1 @> F>> \C\P^1 \\ 
@V \theta  VV @VV \theta   V\\ 
  \C\P^1 @> A>>  \C\P^1. 
\end{CD} 
\ee
commutes. \qed
\et

\subsection{Invariant curves} 
Let $A_1,A_2$ be rational functions, $(A_1,A_2)$ the map given by formulas \eqref{ss0},  \eqref{ss}, and   $C$ an irreducible algebraic curve in $( \C\P^1)^2$.  We say that $C$ is $(A_1,A_2)$-{\it invariant} if  $(A_1,A_2)(  C)=   C.$ We recall that a {\it desingularization} of $  C$ is  a compact Riemann surface  $\widetilde C$
 together with a map $\pi:\widetilde C\rightarrow C$, which is a biholomorphic except for finitely many points. 

The simplest $(A_1,A_2)$-invariant curves are vertical lines $x=a$, where $a$ is a fixed point of $A_1$, and 
horizontal lines $y=b$, where $b$ is a fixed point of $A_2$. Other invariant curves are described as follows (see \cite{inv}, Theorem 4.1).

\bt \l{in1}
Let $A_1,A_2$ be rational functions of degree at least two, and $  C$ an irreducible $(A_1,A_2)$-invariant curve that is not a vertical or horizontal line.
Then the desingularization $\t{  C}$ of $  C$ has genus zero or one, and 
 there exist non-constant holomorphic maps $X_1,X_2:\t{  C}\rightarrow  \C\P^1$ and   $B:\t{  C}\rightarrow \t{  C}$ such that 
 the diagram 
\be \l{du}
\begin{CD} 
(\t{  C})^2 @>(B,B)>> (\t{  C})^2 \\ 
@V (X_1,X_2)  VV @VV  (X_1,X_2) V\\ 
 ( \C\P^1)^2 @>(A_1,A_2)>> ( \C\P^1)^2 
\end{CD} 
\ee
commutes and the map $t\rightarrow (X_1(t),X_2(t))$ is a generically one-to-one parametrization of $  C.$   Finally, unless both $A_1,$ $A_2$ are Latt\`es maps, $\t{  C}$ has genus zero. \qed
\et

For a general description of $(A_1,A_2)$-invariant curves we refer the reader to the paper \cite{inv}. Below we  need only the following description of invariant curves in case $A_1=A_2$  (see \cite{inv}, Theorem 1.2). 

\bt \l{1+} Let $A$ be a rational function of degree at least two that is not a generalized Latt\`es map,  and 
$  C$ an irreducible algebraic curve  in $( \C\P^1)^2$  that is not a vertical or horizontal line. Then $  C$ is $(A,A)$-invariant  if and only if 
 there exist  rational functions  $U_1,$ $U_2,$ $V_1,$ $V_2$    commuting with  $A$ such that 
the equalities 
\be \l{en11} U_1\circ V_1=U_2\circ V_2=A^{\circ d},\ee
 \be\l{en21} V_1\circ U_1=V_2\circ U_2=A^{\circ d}\ee hold for some  $d\geq 0$ and   
the map $t\rightarrow (U_1(t),U_2(t))$ is a parametrization of  $  C$. \qed
\et

Notice that  in general the pa\-ra\-met\-rization  $t\rightarrow (U_1(t),U_2(t))$ provided by Theorem \ref{1+}   is not generically one-to-one. 

\section{Algebraic dependencies between Poincar\'e  functions}\l{sec3}
Our proof of Theorem \ref{m2} is based on the results of Section \ref{s1} and the lemmas below. 

\bl \l{bear} Let $C:\, f(x,y)=0$ be an irreducible algebraic curve that admits a parametrization  $z\rightarrow (\phi_1(z), \phi_2(z))$ by functions  meromorphic on $\C.$  
Then the desingularization $\t{C}$ of $C$ has genus zero or one and there exist meromorphic functions $\phi:\, \C \rightarrow \t{C}$ and $\t \phi_1:\, \t{C}\rightarrow \C\P^1$, $\t \phi_2:\, \t{C}\rightarrow \C\P^1$ such that 
 $$\phi_1=\t \phi_1\circ \phi,\ \  \ \ \  \phi_2=\t \phi_2\circ \phi,$$ and the map 
$z\rightarrow (\t\phi_1(z),\t\phi_2(z))$ from 
$\t{C}$ to $C$ is generically one-to-one. 
\el
\pr The lemma follows  from the Picard theorem
 (see \cite{bn},  Theorem 1 and Theorem 2). \qed

\bl \l{l30}  Let $A$ be a non-special rational function of degree at least two, and  $z_0$ its fixed point with the 
multiplier $\lambda.$ Assume that $W$ is  a rational function of degree at least two commuting with $A$ such that  $z_0$ is  a  fixed point of $W$ with the 
multiplier $\mu.$ Then there exist positive integers $l$ and $k$ such that $\mu^l=\lambda^k.$ 
\el
\pr By the theorem of Ritt, there exist positive integers $l$ and $k$ such that
$W^{\circ l}=A^{\circ k}$, and differentiating this equality at $z_0$ we conclude that  $\mu^l=\lambda^k.$ \qed

\bl \l{xorr} Let $A$, $B$ be rational functions   of degree at least two, and  $X$  a non-constant rational function 
such that the diagram 
\be \l{as}
\begin{CD} 
\C\P^1 @> B>> \C\P^1 \\
@VV X V @VV  X V\\ 
\f \C\P^1 @>A>> \f\C\P^1\
\end{CD}
\ee
commutes. Assume that $z_0$ is a fixed point of $B$ with the multiplier $\lambda_0$. 
Then  $z_1=X(z_0)$ is  a fixed point $z_1$ of $A$ with the multiplier \be \l{eh} \lambda_1=\lambda_0^{\ord_{z_0} X}.\ee In particular, $z_0$ is a repelling fixed point of $B$
if and only if $z_1$ is a repelling fixed point of $A$. Furthermore, if $z_0$ is repelling and  
$\f P_{B,z_0,\lambda}$ is a Poincar\'e function, then  
the equality  
\be \l{sob} {\f P}_{A,z_1,\lambda_1}(z^{\ord_{z_0} X})=X\circ {\f P}_{B,z_0,\lambda_0} \ee
 holds for some Poincar\'e function ${\f P}_{A,z_1,\lambda_1}$. 
\el 
\pr It is clear that $z_1$ is a fixed point of $A$, and a local calculation shows that equality \eqref{eh} holds. Thus, $z_1$ is a repelling fixed point of $A$ if and only if $z_0$ is   a repelling fixed point  of $B$.

The rest of the proof is obtained by a modification of the proof of the uniqueness of a  Poincar\'e function (see e.g. \cite{mil}). Namely, considering the function 
$$G={\f P}_{A,z_1,\lambda_1}^{-1}\circ X\circ {\f P}_{B,z_0,\lambda_0}$$ 
holomorphic in a neighborhood of zero and satisfying $G(0)=0,$ we see 
that  
\begin{multline} G(\lambda_0z)={\f P}_{A,z_1,\lambda_1}^{-1}\circ X\circ B\circ {\f P}_{B,z_0,\lambda_0}=
{\f P}_{A,z_1,\lambda_1}^{-1}\circ A\circ X\circ {\f P}_{B,z_0,\lambda_0}=\\
=\lambda_1\circ {\f P}_{A,z_1,\lambda_1}^{-1}\circ X\circ {\f P}_{B,z_0,\lambda_0}=
\lambda_0^{\ord_{z_0} X}G(z).
\end{multline} 
Comparing now coefficients of the Taylor expansions in the left and the right parts of this equality and taking into account that $\lambda_0$ is not a root of unity, we conclude  
that $G=z^{\ord_{z_0} X},$ implying \eqref{sob}. \qed

\bl \l{l3}
 Let $A$ be a rational function of degree at least two,  $z_0$ its  repelling fixed point with the multiplier $\lambda,$ and 
$\f P_{A,z_0,\lambda}$ a Poincar\'e function. 
 Assume that $C: f(x,y)=0$ is an irreducible algebraic curve, and $d_1,d_2$ are positive integers such that 
the equality 
\be \l{x0} f\left(\f P_{A,z_0,\lambda_0}(z^{d_1}),\f P_{A,z_0,\lambda_0}(z^{d_2})\right)=0\ee  
holds. 
Then $d_1=d_2$, and $C$ is the diagonal $x=y$.  
\el
\pr Since \be \l{sp} z\rightarrow \left(\f P_{A,z_0,\lambda_0}(z^{d_1}),\f P_{A,z_0,\lambda_0}(z^{d_2})\right)\ee is a parametrization of $C$, it is clear that $C$ is not a vertical or horizontal line. Furthermore, substituting $\lambda_0z$ for $z$ into \eqref{x0}, we see that the curve $C$ is $(A^{\circ d_1}, A^{\circ d_2})$-invariant.  
 Therefore, by Theorem \ref{in1}, 
 there exist non-constant holomorphic maps $X_1,X_2:\t{\c C}\rightarrow  \C\P^1$ and   $B:\t{\c C}\rightarrow \t{\c C}$ such that 
 the diagram 
\be \l{du1}
\begin{CD} 
(\t{\c C})^2 @>(B,B)>> (\t{\c C})^2 \\ 
@V (X_1,X_2)  VV @VV  (X_1,X_2) V\\ 
 ( \C\P^1)^2 @>(A^{\circ d_1},A^{\circ d_2})>> ( \C\P^1)^2 
\end{CD} 
\ee
 commutes.  Thus, 
$$\deg A^{\circ d_1}=\deg A^{\circ d_2}=\deg B,$$ and hence $d_1=d_2.$  Since the parametrization of $C$ has the form \eqref{sp}, this implies  
that $C$ is the diagonal. 
\qed 

\bc \l{c3}  Let $A_1$, $A_2$ be  rational functions of degree at least two,  $z_1$, $z_2$ their repelling fixed points with multipliers $\lambda_1,$ $\lambda_2$, and 
$\f P_{A_1,z_1,\lambda_1}$, $\f P_{A_2,z_2,\lambda_2}$  Poincar\'e functions. 
 Assume that $C: f(x,y)=0$ is an irreducible algebraic curve  and $d_1,d_2,\t d_1,\t d_2$ are positive integers such that $\GCD(d_1,d_2)=1$ and 
the equalities 
\be \l{x1} f\left(\f P_{A_1,z_1,\lambda_1}(z^{d_1}),\f P_{A_2,z_2,\lambda_2}(z^{d_2})\right)=0,\ee  
\be \l{x2} f\left(\f P_{A_1,z_1,\lambda_1}(z^{\t d_1}),\f P_{A_2,z_2,\lambda_2}(z^{\t d_2})\right)=0\ee 
hold. 
Then there exists a positive integer $k$ such that the equalities 
\be \l{suth} \t d_1=kd_1, \ \ \ \ \t d_2=kd_2\ee hold.    
\ec
\pr It is clear that equalities \eqref{x1}, \eqref{x2} imply the 
equalities 
\be \l{xx1} f\left(\f P_{A_1,z_1,\lambda_1}(z^{d_1\t d_1}),\f P_{A_2,z_2,\lambda_2}(z^{d_2\t d_1})\right)=0\ee and 
\be \l{xx2} f\left(\f P_{A_1,z_1,\lambda_1}(z^{d_1\t d_1}),\f P_{A_2,z_2,\lambda_2}(z^{d_1\t d_2})\right)=0. \ee
Eliminating now from these equalities $\f P_{A_1,z_1,\lambda_1}(z^{d_1\t d_1})$, we conclude that the functions $\f P_{A_2,z_2,\lambda_2}(z^{d_2\t d_1})$ and $\f P_{A_2,z_2,\lambda_2}(z^{d_1\t d_2})$ are algebraically dependent. Therefore, $\t d_1d_2=d_1\t d_2$ by Lemma \ref{l3}, implying \eqref{suth}. \qed

\vskip 0.2cm
\noindent{\it Proof of Theorem \ref{m2}.} 
Let $C: f(x,y)=0$ be an irreducible algebraic curve with a generically one-to-one parametrization by rational functions 
$z\rightarrow (X_1(z),X_2(z))$, and $d_1,$ $d_2$ coprime positive integers. 
Assume that diagram \eqref{hol} commutes for some rational function $B$ with a repelling fixed point $z_0$  and equalities \eqref{io}, \eqref{oi} hold. Then denoting the multiplier of $z_0$ by $\lambda$ and using   Lemma \ref{xorr}, we see that 
\be \l{korova} \lambda_1^{l_1}=\lambda^{\ord_{z_0}X_1}, \ \ \ \ \lambda_2^{l_2}=\lambda^{\ord_{z_0}X_2},\ee 
and 
\begin{multline} \l{zaq0} 
0=f(X_1,X_2)=f(X_1\circ \f P_{B,z,\lambda},X_2\circ \f P_{B,z,\lambda})=\\ 
=f\Big(\f P_{A_1^{\circ l_1},z_1,\lambda_1^{l_1}}\left(z^{\ord_{z_0}X_1}\right),\f P_{A_2^{\circ l_2},z_2,\lambda_2^{l_2}}\left(z^{\ord_{z_0}X_2}\right)\Big).
\end{multline}
Since 
$$
\f P_{A_1^{\circ l_1},z_1,\lambda_1^{l_1}}(z)= \f P_{A_1,z_1,\lambda_1}(z), \ \ \ \ \f P_{A_2^{\circ l_2},z_2,\lambda_2^{l_2}}(z)= \f P_{A_2,z_2,\lambda_2}(z),
$$
this implies  that 
\be \l{zaq} f\Big(\f P_{A_1,z_1,\lambda_1}\left(z^{\ord_{z_0}X_1}\right),\f P_{A_2,z_2,\lambda_2}\left(z^{\ord_{z_0}X_2}\right)\Big)=0.\ee
Finally, \eqref{oi} implies that if \eqref{zaq} holds, then \eqref{xru} also holds. 
This proves the ``if'' part of the theorem.

To prove the ``only if'' part, it is enough to show that equality \eqref{xru} implies that there exist positive integers $r_1,r_2$ such that 
\be \l{tc} \lambda_1^{r_1}= \lambda_2^{r_2}=\lambda.\ee Indeed, in this case 
 substituting $\lambda z$ for $z$ into \eqref{xru} we obtain the equality 
$$f\Big(A_1^{\circ d_1r_1}\circ \f P_{A_1,z_1,\lambda_1}(z^{d_1}), A_2^{\circ d_2r_2}\circ\f P_{A_2,z_2,\lambda_2}(z^{d_2})\Big)=0.$$ Therefore, for $$l_1=d_1r_1,\ \ \ l_2=d_2r_2,$$ 
 the curve $C$ is $(A_1^{\circ l_1}, A_2^{\circ l_2})$-invariant,  implying by Theorem \ref{in1} that $C$ has genus zero and there exist rational functions $X_1,X_2$ and $B$ such that diagram \eqref{hol} commutes and  the map $z\rightarrow (X_1(z),X_2(z))$ 
is a generically one-to-one parametrization of $C$. 
It follows now from Lemma \ref{bear} that there exists 
a meromorphic function $\phi$ such that  the equalities 
\be \l{tgv} \f P_{A_1,z_1,\lambda_1}(z^{d_1})=X_1\circ \phi(z), \ \ \ \ \f P_{A_2,z_2,\lambda_2}(z^{d_2})=X_2\circ \phi(z).\ee 
hold. 
Thus, 
$$z_1=\f P_{A_1,z_1,\lambda_1}(0)=X_1\circ \phi(0), \ \ \ \ z_2=\f P_{A_2,z_2,\lambda_2}(0)=X_2\circ \phi(0)
,$$ implying that equalities \eqref{io} hold for the point $z_0=\phi(0).$ 

Further, since $z_1$ and $z_2$ are fixed points of $A_1$ and $A_2$, the point $z_0$ is a preperiodic point 
of $B$. Thus, changing in \eqref{hol} the functions $B$ and $A_1^{\circ l_1},$ $A_2^{\circ l_2}$ to some of their iterates, and the point $z_0$ to some point in its $B$-orbit,  we may assume that $z_0$ is a fixed point of $B$. Moreover, $z_0$ is repelling by Lemma \ref{xorr}. Let us recall now that, by what is proved above, \eqref{hol} and \eqref{io} imply \eqref{zaq}. Thus,  equalities \eqref{xru} and \eqref{zaq} hold simultaneously and hence equalities \eqref{oi}  hold   
by Corollary \ref{c3}.

Let us show now that \eqref{xru} implies \eqref{tc}. Assume first that $A_1$ and $A_2$  are not generalized Latt\`es maps. Substituting $\lambda_2z$ for $z$ into  equality \eqref{xru} we obtain  
the equality 
\begin{multline}
f\Big(\f P_{A_1,z_1,\lambda_1}\circ (\lambda_2 z)^{d_1}, \f P_{A_2,z_2,\lambda_2}\circ (\lambda_2z)^{d_2}\Big)=\\
 =f\Big(\f P_{A_1,z_1,\lambda_1}\circ (\lambda_2 z)^{d_1}, A_2^{\circ d_2}\circ \f P_{A_2,z_2,\lambda_2}\circ z^{d_2}\Big)=
0,
\end{multline}
implying that the functions $\f P_{A_1,z_1,\lambda_1}\circ (\lambda_2 z)^{d_1}$ and $\f P_{A_2,z_2,\lambda_2}\circ z^{d_2}$ satisfy the equality 
\be \l{krot} g\Big(\f P_{A_1,z_1,\lambda_1}\circ (\lambda_2 z)^{d_1}, \f P_{A_2,z_2,\lambda_2}\circ z^{d_2}\Big)=
0,\ee where $g(x,y)=f(x,A_2^{\circ d_2}(y)).$  Eliminating now from \eqref{xru} and \eqref{krot} the function
 $\f P_{A_2,z_2,\lambda_2}\circ z^{d_2}$, we conclude that 
the functions $\f P_{A_1,z_1,\lambda_1}\circ z^{d_1}$ and $\f P_{A_1,z_1,\lambda_1}\circ (\lambda_2 z)^{d_1}$ are  algebraically dependent. In turn, this implies that the functions $\f P_{A_1,z_1,\lambda_1}(z)$ and  $\f P_{A_1,z_1,\lambda_1}(\lambda_2^{d_1} z)$  also are  algebraically dependent.

Let $\widetilde C: \widetilde f(x,y)=0$ be a curve such that 
$$\widetilde f\Big(\f P_{A_1,z_1,\lambda_1}(z), \f P_{A_1,z_1,\lambda_1}(\lambda_2^{d_1}z)\Big)=0.$$ Then  substituting $\lambda_1 z$ for $z$ we see that $\widetilde f$ is $(A_1,A_1)$-invariant. Therefore, by Theorem \ref{1+}, there exist rational function $V_1$ and $V_2$ commuting with  $A_1$ such that $\widetilde C$ is a component of the curve $$V_1(x)-V_2(y)=0,$$ implying that  the equality 
\be \l{th} V_1\circ \f P_{A_1,z_1,\lambda_1}(z)=V_2\circ  \f P_{A_1,z_1,\lambda_1}(\lambda_2^{d_1} z)\ee holds. Furthermore, it follows from the Ritt theorem  that there exist positive integers $s_1,s_2$, and $s$ such that \be \l{lu} V_1^{\circ s_1}=V_2^{\circ s_2}=A_1^{\circ s}.\ee
Since \eqref{th} implies that for every $l\geq 1$ the equality 
$$V_1^{\circ l}\circ V_1\circ \f P_{A_1,z_1,\lambda_1}(z)=V_1^{\circ l}\circ V_2\circ  \f P_{A_1,z_1,\lambda_1}(\lambda_2^{d_1} z)$$ holds, 
setting 
$$W_1=V_1^{\circ s_1}, \ \ \ W_2=V_1^{\circ (s_1-1)}\circ V_2,$$ we see that 
 $W_1$ and $W_2$ also commute with $A_1$ and satisfy 
\be \l{mor} W_1\circ \f P_{A_1,z_1,\lambda_1}(z)=W_2\circ  \f P_{A_1,z_1,\lambda_1}(\lambda_2^{d_1} z).\ee 
In addition, $z_1$ is a fixed point 
of $W_1$ by \eqref{lu}. Finally, since equality \eqref{mor} implies the equality $$W_1(z_1)=W_2(z_1),$$ the point $z_1$ is also a   fixed point of $W_2$.

Differentiating equality \eqref{mor} at zero, we see that  the multipliers
$$\mu_1=W_1'(z_1), \ \ \ \mu_2=W_2'(z_1)$$ satisfy the equality 
\be \l{if} \mu_1=\mu_2\lambda_2^{d_1}.\ee On the other hand,  Lemma \ref{l30} yields that   there exist positive integers $k_1,$ $k_2,$ and $k$ such that 
\be \l{fi} \mu_1^{k_1}=\mu_2^{k_2}=\lambda_1^k.\ee 
It follows now from \eqref{if} and \eqref{fi} that 
$$\lambda_1^{kk_2}=\mu_1^{k_1k_2}=\mu_2^{k_1k_2}\lambda_2^{d_1k_1k_2}=\lambda_1^{kk_1}\lambda_2^{d_1k_1k_2},$$ implying 
that 
\be \l{ei} \lambda_1^{k(k_2-k_1)}=\lambda_2^{d_1k_1k_2}.\ee
Moreover, since $\vert \lambda_1 \vert >1,$ $ \vert \lambda_2 \vert >1,$  the number $k_2- k_1$ is positive.  This proves the implication 
\eqref{xru}$\Rightarrow$\eqref{tc} in case $A_1$ and $A_2$  are not generalized Latt\`es maps.

Assume now that    $A_1$, $A_2$ are arbitrary non-special rational functions. Then, by Theorem \ref{33},  there exist rational functions $F_1$, $F_2,$ $\theta_1,$ $\theta_2$ such that the diagrams 
\be 
\begin{CD} 
\C @>F_1>> \C \\
@VV \theta_1 V @VV \theta_1 V\\ 
\C\P^1 @>A_1>> \C\P^1\ ,
\end{CD}
\ \ \ \ \ \ \ \ 
\begin{CD} 
\C @>F_2>> \C \\
@VV \theta_2 V @VV \theta_2 V\\ 
\C\P^1 @>A_2>> \C\P^1\
\end{CD}
\ee
commute, and $F_1$, $F_2$ are 
not generalized Latt\`es maps. Further, since all the points in the preimage $\theta_{A_i}^{-1}\{z_i\},$ $i=1,2$, are $F_i$-preperiodic, there exist  a positive integer $N$ and fixed points $z_1',$ $z_2'$ of $F_1^{\circ N},$ $F_2^{\circ N}$ such that the diagrams 
\be 
\begin{CD} 
\C @>F_1^{\circ N}>> \C \\
@VV \theta_1 V @VV \theta_1 V\\ 
\C\P^1 @>A_1^{\circ N}>> \C\P^1\ ,
\end{CD}
\ \ \ \ \ \ \ \ 
\begin{CD} 
\C @>F_2^{\circ N}>> \C \\
@VV \theta_2 V @VV \theta_2 V\\ 
\C\P^1 @>A_2^{\circ N}>> \C\P^1\
\end{CD}
\ee
commute, and the equalities  
$$\theta_1(z_1')=z_1, \ \ \ \ \ \theta_1(z_2')=z_2$$ hold. Moreover, if $\mu_i$ is 
the multiplier of $F_i^{\circ N}$ at $z_i'$, $i=1,2,$ then, by Lemma \ref{xorr}, the equalities 
\be \l{ws} \mu_1^{\ord_{z_1'}\theta_1}=\lambda_1^{N}, \ \ \ \mu_2^{\ord_{z_2'}\theta_2}=\lambda_2^{N},\ee   
\be \l{se1} \f P_{A_1^{\circ N},z_1,\lambda_1^N}(z^{\ord_{z_1'}\theta_1})=\theta_1\circ \f P_{F_1^{\circ N},z_1',\mu_1}(z),\ee \be \l{se2} \f P_{A_2^{\circ N},z_2,\lambda_2^N}(z^{\ord_{z_2'}\theta_2})=\theta_2\circ \f P_{F_2^{\circ N},z_2',\mu_2}(z)\ee hold. 
  
Setting 
$$f_1=\ord_{z_1'}\theta_1, \ \ \ \ 	f_2=\ord_{z_2'}\theta_2, \ \ \ \ f=f_1f_2,$$ 
and substituting  $z^{d_1f_2}$ and $z^{d_2f_1}$ for $z$  into 
equalities \eqref{se1} and \eqref{se2},  we obtain that     
$$ \f P_{A_1,z_1,\lambda_1}(z^{d_1f})=\f P_{A_1^{\circ N},z_1,\lambda_1^N}(z^{d_1f})=\theta_1\circ \f P_{F_1^{\circ N},z_1',\mu_1}(z^{d_1f_2}),$$
$$ \f P_{A_2,z_2,\lambda_2}(z^{d_2f})= \f P_{A_2^{\circ N},z_2,\lambda_2^N}(z^{d_2f})=\theta_2\circ \f P_{F_2^{\circ N},z_2',\mu_2}(z^{d_2f_1}).$$ 
Thus, equality \eqref{xru} implies that the functions $\f P_{F_1^{\circ N},z_1',\mu_1}(z^{d_1f_2})$ and $\f P_{F_2^{\circ N},z_2',\mu_2}(z^{d_2f_1})$ satisfy the equality 
$$\t f\left(\f P_{F_1^{\circ N},z_1',\mu_1}(z^{d_1f_2}),\f P_{F_2^{\circ N},z_2',\mu_2}(z^{d_2f_1})\right)=0,$$ where $$\t f(x,y)=f\left(\theta_1(x),\theta_2(y)\right).$$
Since $F_1^{\circ N}, F_2^{\circ N}$ are not generalized Latt\`es maps, by what is proved above there   
exist positive integers $p_1,p_2$ such that 
$\mu_1^{p_1}= \mu_2^{p_2},$ implying by \eqref{ws} that 
$$\lambda_1^{p_1f_2N}=\mu_1^{p_1f_1f_2}=\mu_2^{p_2f_1f_2}=\lambda_2^{p_2f_1N}.$$
Thus, equality 
\eqref{tc} holds for the integers
$$r_1=p_1f_2N, \ \ \ \ r_2=p_2f_1N. \eqno{\Box}$$


\vskip 0.2cm
\noindent{\it Proof of Corollary \ref{cor1}.} 
If $\f P_{A_1,z_1,\lambda_1}$, $\f P_{A_2,z_2,\lambda_2}$ are algebraically dependent, then 
it follows from the commutativity of  diagram \eqref{hol} that 
$$(\deg A_1)^{l_1}=(\deg A_2)^{l_2}=\deg B,$$ implying that $n_1^{l_1}=n_2^{l_2}.$ 
Furthermore, it follows from equalities \eqref{korova} that 
 $$\lambda_1^{l_1\ord_{z_0}X_2}=\lambda_2^{l_2\ord_{z_0}X_1}. \eqno{\Box} $$

\vskip 0.2cm

The following result shows that if $A_1$ and $A_2$ are not generalized Latt\`es maps, then dependencies \eqref{xru} actually reduce to dependencies \eqref{xre}.

\bt \l{m1} Let $A_1$, $A_2$ be  rational functions of degree at least two   that are not generalized Latt\`es maps,  $z_1$, $z_2$ their repelling fixed points with multipliers $\lambda_1,$ $\lambda_2$,  
and $\f P_{A_1,z_1,\lambda_1}$, $\f P_{A_2,z_2,\lambda_2}$ Poincar\'e functions. 
Assume that  $C: f(x,y)=0$ is an irreducible algebraic curve, and $d_1,$ $d_2$ are coprime positive integers such that the equality 
\be \l{xrus} f\left(\f P_{A_1,z_1,\lambda_1}(z^{d_1}),\f P_{A_2,z_2,\lambda_2}(z^{d_2})\right)=0\ee holds. Then  $d_1=d_2=1$  and $C$ has genus zero. 
Furthermore, if  $C: f(x,y)=0$ is an irreducible curve of genus zero with a generically one-to-one parametrization by rational functions 
$z\rightarrow (X_1(z),X_2(z))$, then the equality 
$$
f\left(\f P_{A_1,z_1,\lambda_1}(z),\f P_{A_2,z_2,\lambda_2}(z)\right)=0$$
 holds for some  Poincar\'e functions $\f P_{A_1,z_1,\lambda_1}$, $\f P_{A_2,z_2,\lambda_2}$ if and only if  
 there exist positive integers $l_1,l_2$ and a rational function $B$ with a repelling fixed point  $z_0$  such that the diagram 
\be \l{hols} 
\begin{CD} 
 (\C\P^1)^2 @> (B,B) >> (\C\P^1)^2
 \\ 
@V (X_1,X_2) VV @VV  (X_1,X_2)   V\\ 
(\C\P^1)^2 @> (A_1^{\circ l_1}, A_2^{\circ l_2}) >> (\C\P^1)^2, 
\end{CD} 
\ee
commutes, and the equalities 
\be \l{ios} X_1(z_0)=z_1, \ \ \ \ X_2(z_0)=z_2,\ee  
\be \l{kon} X_1'(z_0)\neq 0, \ \ \  \ X_2'(z_0)\neq 0\ee
hold. 
\et 
\pr The proof if obtained by a modification of the proof of Theorem \ref{m2}, taking into account that if  $A_1,$ $A_2$ are not generalized Latt\`es maps, then it follows from the commutativity of diagram \eqref{hol} by Proposition \ref{1} that there exist rational functions $Y_1$ and $Y_2$ such that the equalities 
$$  Y_1  \circ X_1  =B^{\circ d_1} \ \ \ \ \ \ \ Y_2 \circ X_2 =B ^{\circ d_2}
$$ hold for some $d_1,d_2\geq 0$.
Therefore, for any repelling fixed point  $z_0$ of $B$ the inequalities \eqref{kon} hold by the chain rule.
Thus, $d_1=d_2=1$ by \eqref{oi}. 
 \qed  
\vskip 0.2cm
Notice that unlike the case of B\"ottcher functions, algebraic dependencies \eqref{xre} of degree greater than one between Poincar\'e functions do exist. The simplest of them are graphs constructed as follows. Let us take any two rational functions $U$ and $V$, and set 
\be \l{edc} A_1=U\circ V, \ \ \ \ \ A_2=V\circ U.\ee Then 
the diagram 
\be 
\begin{CD} 
 \C\P^1 @> A_1 >> \C\P^1
 \\ 
@V {V} VV @VV  {V}  V\\ 
\C\P^1 @> A_2 >> \C\P^1 
\end{CD} 
\ee 
obviously commutes. Moreover, if $z_0$ is a repelling fixed point of $A_1$, then the point $z_1=V(z_0)$ is a repelling 
fixed point of $A_2$ by Lemma \ref{xorr}. 
Finally, the first equality in \eqref{edc} implies that $V'(z_1)\neq 0.$ 
Therefore,  
$$\f P_{A_2,z_2,\lambda_2}=V\circ \f P_{A_1,z_1,\lambda_1}, $$
by Lemma \ref{xorr}.

Notice also that the equality $d_1=d_2=1$ provided by Theorem \ref{m1} does not hold for arbitrary
 non-special $A_1,$ $A_2$. For example, let $A$ be any rational function  of the form $A=zR^d(z)$, where $R\in \C(z)$ and $d>1.$ Then one can easily check that
 $A:\f O\rightarrow \f O$, where $\f O$ is defined by the equalities 
$$\nu(0)=d, \ \ \ \nu(\infty)=d,$$ is a minimal holomorphic map between orbifolds. Thus, $A$ is a generalized Latt\`es map. Furthermore, 
the diagram 
$$
\begin{CD}
\C\P^1 @>zR(z^d)>> \C\P^1 \\
@VVz^d V @VVz^d V\\ 
\C\P^1 @> zR^d(z)>> \ \ \C \P^1.
\end{CD}
$$
obviously commutes. Choosing now $R$ in such a way that zero  is a repelling fixed point of $zR(z^d)$ and denoting  by $\lambda$ the multiplier of $zR^d(z)$ at zero, we obtain 
by Lemma \ref{xorr} that  
 $$
 {\f P}_{zR^d(z),0,\lambda^d}(z^{d})=z^d\circ {\f P}_{zR(z^d),0,\lambda}(z).
$$
 Thus, $ {\f P}_{zR^d(z),0,\lambda^d}(z^{d})$ and ${\f P}_{zR(z^d),0,\lambda}(z)$ are algebraically dependent. 

\section{Algebraic dependencies between B\"ottcher functions} 
\subsection{Polynomial semiconjugacies and invariant curves} 

If $A_1$, $A_2$ are non-special  {\it polynomials} of degree at least two, then  any  irreducible $(A_1,A_2)$-invariant curve $C$  that is not a vertical or horizontal line has genus zero  and allows for a generically one-to-one  parametrization by {\it polynomials}  $X_1,X_2$ such that 
 the diagram 
\be \l{dur}
\begin{CD} 
  ( \C\P^1)^2 @>(B,B)>>  ( \C\P^1)^2 \\ 
@V (X_1,X_2)  VV @VV  (X_1,X_2) V\\ 
 ( \C\P^1)^2 @>(A_1,A_2)>> ( \C\P^1)^2 
\end{CD} 
\ee
commutes for some {\it polynomial} $B$ (see Proposition 2.34
of \cite{ms} or Section 4.3 of \cite{pj}).

For fixed polynomials $A$, $B$ of degree at least two, we denote  by $\f E(A,B)$ the set (possibly empty) consisting of polynomials $X$ of degree at least two such that diagram \eqref{ssee} 
 commutes. The following result was proved in the paper \cite{pj} as a corollary of results of the paper \cite{p1}.

\bt \l{univ} Let $A$ and $B$ be fixed  non-special polynomials of degree at least two
such that the set $\f E(A,B)$ is non-empty, and let $X_0$ be an element of  $\f E(A,B)$ of the minimum possible degree. Then a polynomial $X$ belongs to 
$\f E(A,B)$ if and only if $X= \t A\circ X_0$ for some polynomial  $\t A$ commuting with $A.$ \qed 
\et 

Notice that applying Theorem \ref{univ} for $B=A$ one can  reprove the classification of commuting polynomials  and, more generally, of commutative semigroups of $\C[z]$ obtained  in the papers \cite{r3}, \cite{r}, \cite{wora}
(see \cite{amen}, Section 7.1,  for more detail).  
On the other hand, applying Theorem \ref{univ} to system \eqref{dur} with $A_1=A_2=A$, we see that $X_1,$ $X_2$ cannot provide a generically one-to-one  parametrization of $C$, unless one of the polynomials $X_1$, $X_2$ has degree one. Moreover if, say, $X_1$ has degree one, then  without loss of generality we may assume that $X_1=z,$  implying that $B=A$ and $X_2$ commutes with $A$. Thus, we obtain the following result obtained  by Medve\-dev and Scanlon in the paper \cite{ms}. 

\bt \l{med} Let $A$ be a non-special polynomial of degree at least two, and 
$C$  an irreducible algebraic curve that is not a vertical or horizontal line. Then $C$ is $(A,A)$-invariant if and only if $C$ has the form $x= P(y)$ or $y=P(x)$, where  $P$ is  a polynomial commuting with $A$. \qed
\et

Finally, yet another corollary of Theorem \ref{univ} is the following result, which complements the classification of $(A_1,A_2)$-invariant curves obtained in \cite{ms} 
(see \cite{pj}, Theorem 1.4). 
\bt \l{new} Let $A_1$, $A_2$ be non-special polynomials of degree at least two, and $C$ a curve. Then $C$ is an irreducible $(A_1,A_2)$-invariant curve if and only if 
 $C$ has the form $Y_1(x)-Y_2(y)=0$, where $Y_1,Y_2$  are  polynomials of coprime degrees satisfying  the equations 
\be \l{krys} T\circ Y_1=Y_1\circ A_1, \ \ \ \ T\circ Y_2=Y_2\circ A_2\ee for some polynomial $T.$ \qed
\et

\subsection{Proof of Theorem \ref{m3}}
As in the case of  Poincar\'e functions, we do not assume that considered 
B\"ottcher functions are normalized. Thus, 
the notation $\f B_P$ is used to denote  
{\it some}  function satisfying  conditions \eqref{boet1}, \eqref{boet2}. 

To prove Theorem \ref{m3} we need the following two lemmas. 

\bl \l{xerr} Let $A$, $B$ be polynomials     of degree at least two, and  $X$  a non-constant polynomial  such that the diagram 
\be \l{ass}
\begin{CD} 
\C\P^1 @> B>> \C\P^1 \\
@VV X V @VV  X V\\ 
\f \C\P^1 @>A>> \f\C\P^1\
\end{CD}
\ee
commutes. Assume that $\f B_{B}$ is a B\"ottcher function. Then   
\be \l{sobb} X\circ {\f B}_B(z)= {\f B}_{A}(z^{\deg  X})\ee
for some B\"ottcher function ${\f B}_{A}$.  
\el 
\pr The lemma follows from  Lemma 2.1 of \cite{ng}. \qed

\bl \l{l32}  Let $A$ be a polynomial of degree $n\geq 2$,  and 
$\f B_{A}$ a B\"ottcher function. 
 Assume that $C: f(x,y)=0$ is an irreducible algebraic curve and $d_1,d_2$ are positive integers such that $d_1\leq d_2$ and 
the equality 
\be \l{x01} f\left(\f B_{A}(z^{d_1}),\f B_{A}(z^{d_2})\right)=0\ee  
holds. Then $C$ is a graph \be \l{gra} P(x)-y=0,\ee where $P$ is a polynomial commuting with $A$, 
and the equality \be \l{sle} d_1\,\deg P=d_2 \ee holds.
\el
\pr Substituting $z^n$ for $z$ in \eqref{x01}, we see that the curve $C$ is $(A, A)$-invariant.
Therefore, by Theorem \ref{med},  $C$ is a graph of the form $x= P(y)$ or $y=P(x)$, where  $P$ is  a polynomial commuting with $A$. Taking into account that $d_1\leq d_2$, this implies that  \eqref{gra} and  
\eqref{sle} hold. \qed

\bc \l{c32}  Let $A_1$, $A_2$ be  polynomials of degree at least two, and 
$\f B_{A_1}$, $\f B_{A_2}$ B\"ottcher functions. 
 Assume that $C: f(x,y)=0$ is an irreducible algebraic curve of genus zero and $d_1,d_2,\t d_1,\t d_2$ are positive integers such that $\GCD(d_1,d_2)=1$ and 
the equalities 
\be \l{x11} f\left(\f B_{A_1}(z^{d_1}),\f B_{A_2}(z^{d_2})\right)=0,\ee  
\be \l{x21} f\left(\f B_{A_1}(z^{\t d_1}),\f B_{A_2}(z^{\t d_2})\right)=0\ee 
hold. Then there exists a positive integer $k$ such that the equalities 
\be \l{suth+} \t d_1=kd_1, \ \ \ \ \t d_2=kd_2\ee hold.   
\ec
\pr It is clear that equalities \eqref{x11}, \eqref{x21} imply the 
equalities 
\be \l{xx11} f\left(\f B_{A_1}(z^{d_1\t d_1}),\f B_{A_2}(z^{d_2\t d_1})\right)=0\ee and 
\be \l{xx21} f\left(\f B_{A_1}(z^{d_1\t d_1}),\f B_{A_2}(z^{d_1\t d_2})\right)=0, \ee
 and eliminating from these equalities the function $\f B_{A_1}(z^{d_1\t d_1})$, we conclude that 
 the functions $\f B_{A_2}(z^{d_2\t d_1})$ and $\f B_{A_2}(z^{d_1\t d_2})$ are algebraically dependent. Therefore, by Lemma \ref{l32},  
one of these functions is a polynomial in the other. 

Assume, say, that 
\be \l{kjh} 
\f B_{A_2}(z^{d_2\t d_1})=R\circ \f B_{A_2}(z^{d_1\t d_2})
\ee
(the other case is considered similarly). Then substituting the right part of this equality for the
 left part into \eqref{xx11}, we conclude that 
$$ f\left(\f B_{A_1}(z^{d_1\t d_1}),R\circ \f B_{A_2}(z^{d_1\t d_2})\right)=0,$$ implying that 
\be \l{ok}  f\left(\f B_{A_1}(z^{\t d_1}),R\circ \f B_{A_2}(z^{\t d_2})\right)=0.\ee
Let us observe now that equalities \eqref{x21} and \eqref{ok} imply that the curve $f(x,y)=0$ is 
invariant under the map $$(z_1,z_2)\rightarrow (\widehat A_1(z_1),\widehat A_2(z_2))=(z_1,R(z_2)).$$ Since the commutativity of \eqref{dur} implies that $\deg A_1=\deg A_2,$ this yields that    $\deg R=1$. It follows now from \eqref{kjh} that  
$$d_2\t d_1=d_1\t d_2,$$ implying \eqref{suth+}. 
\qed

\vskip 0.2cm
\noindent{\it Proof of Theorem \ref{m3}.} 
To prove the ``if'' part, it is enough to observe that if \eqref{hol2} and \eqref{oi2} hold, then by Lemma \ref{xerr} we have
$$0=f(X_1,X_2)=f(X_1\circ {\f B}_B(z),X_2\circ {\f B}_B(z))=f({\f B}_{A_1}^{\circ l_1}(z^{\deg  X_1}),{\f B}_{A_2}^{\circ l_2}(z^{\deg  X_2}))=$$
$$=f({\f B}_{A_1}(z^{\deg  X_1}),{\f B}_{A_2}(z^{\deg  X_2})).$$

In the other direction, if \eqref{xru2} holds, then setting $n_1=\deg A_1,$ $n_2=\deg A_2$, and substituting $z^{n_2}$ for $z$ into \eqref{xru2} we obtain the equality 
\be \l{elim} f(\f B_{A_1}(z^{d_1n_2}), A_2\circ \f B_{A_2}(z^{d_2}))=0.\ee 
 Eliminating now $\f B_{A_2}(z^{d_2})$ from \eqref{xru2} and \eqref{elim}, we conclude that the  functions $\f B_{A_1}(z^{d_1})$ and  $\f B_{A_1}(z^{d_1n_2})$ are algebraically dependent. Since the corresponding algebraic curve $\t f(x,y)=0$ such that 
\be \l{til} \t f(\f B_{A_1}(z^{d_1}),  \f B_{A_1}(z^{d_1n_2}))=0\ee  
is $(A_1,A_1)$-invariant, it follows from Theorem \ref{med} that  
\be \l{pos} \f B_{A_1}(z^{d_1n_2})=P\circ \f B_{A_1}(z^{d_1})\ee for some 
polynomial $P$ commuting with $A_1.$ Clearly, equality \eqref{pos} implies that $\deg P=n_2.$ On the other hand,  by the Ritt theorem, $P$ and $A_1$ have a common iterate. Therefore,  
there exist positive integers $l_1,$ $l_2$ such $n_1^{l_1}=n_2^{l_2}.$ 

Setting now 
$$n=n_1^{l_1}=n_2^{l_2}$$ and substituting $z^{n}$ for $z$ into \eqref{xru2} we obtain that 
$f(x,y)=0$ is $(A_1^{\circ l_1}, A_2^{\circ l_2})$-invariant, implying that  \eqref{hol2} holds. Moreover, by Theorem \ref{new},  $f(x,y)=0$ has the form 
\be \l{ofs} Y_1(x)-Y_2(y)=0,\ee where $Y_1,Y_2$  are  polynomials of coprime degrees. 
 Since a generically one-to-one parametrization $z\rightarrow (X_1(z),X_2(z))$ of \eqref{ofs} satisfies the conditions 
$$\deg X_1=\deg Y_2,\ \ \ \ \deg X_2=\deg Y_1,$$ we conclude that the degrees  
$$\deg X_1=d_1', \ \ \ \ \deg X_2=d_2'$$ of the functions $X_1$ and $X_2$ in \eqref{hol2} 
satisfy $\GCD(d_1',d_2')=1.$ 
Using now the ``if'' part of the theorem, we 
see that equalities \eqref{xru2} and 
$$ f(\f B_{A_1}(z^{d_1'}), \f B_{A_2}(z^{d_2'}))=0$$ 
hold simultaneously, implying   by Corollary \ref{c32} that equalities $d_1'=d_1,$ $d_2'=d_2,$ and \eqref{oi2} hold. \qed

\end{document}